\documentclass[11pt]{amsart}

\theoremstyle{plain}
\newtheorem{theorem}                 {Theorem}      [section]
\newtheorem{proposition}  [theorem]  {Proposition}

\newtheorem*{theorem*}               {Theorem\;{\bf \ref{th: hyper_S4}}}

\theoremstyle{definition}

\newtheorem{remark}       [theorem]  {Remark}

\numberwithin{equation}{section}

\def \1{\mbox{${\mathbf 1}$}}
\def \r{\mbox{${\mathbb R}$}}
\def \s{\mbox{${\mathbb S}$}}

\def \e{\mbox{${\mathbb E}$}}
\def \h{\mbox{${\mathbb H}$}}

\def \rrm{\mbox{${\scriptstyle \frac{1}{\sqrt 2}}$}}

\DeclareMathOperator{\trace}{trace}
\DeclareMathOperator{\grad}{grad}

\DeclareMathOperator{\Div}{div}

\begin{document}

\title{Biharmonic hypersurfaces in $4$-dimensional space forms}
\author{A. Balmu\c s}
\address{Universit\`a degli Studi di Cagliari\\
Dipartimento di Matematica\\
\newline
Via Ospedale 72\\
09124 Cagliari, ITALIA} \email{balmus@unica.it, montaldo@unica.it}
\author{S. Montaldo}
\author{C. Oniciuc}
\address{Faculty of Mathematics, ``Al.I.~Cuza'' University of Iasi\\
\newline
Bd. Carol I Nr. 11 \\
700506 Iasi, ROMANIA} \email{oniciucc@uaic.ro}

\dedicatory{}


\thanks{The first author was supported by a INdAM doctoral fellowship, ITALY.
The second author was supported by PRIN, ITALY. The third author was
supported by Grant CEEX ET 5871/2006, ROMANIA}

\begin{abstract}
We investigate proper biharmonic hypersurfaces with at most three
distinct principal curvatures in space forms. We obtain the full
classification of  proper biharmonic hypersurfaces in
$4$-dimensional space forms.
\end{abstract}


\maketitle

\section{Introduction}

{\it Biharmonic maps} $\varphi:(M,g)\to(N,h)$ between Riemannian
manifolds are critical points of the {\em bienergy} functional
$$
E_2(\varphi)=\frac{1}{2}\int_{M}\,|\tau(\varphi)|^2\,v_g,
$$
where $\tau(\varphi)=\trace\nabla d\varphi$ is the {tension field}
of $\varphi$ that vanishes for harmonic maps. The Euler-Lagrange
equation corresponding to $E_2$ is given by the vanishing of the
{\it bitension field}
$$
\tau_2(\varphi)=-J^{\varphi}(\tau(\varphi))=-\Delta\tau(\varphi)
-\trace \ R^N(d\varphi,\tau(\varphi))d\varphi,
$$
where $J^{\varphi}$ is formally the Jacobi operator of $\varphi$.
The operator $J^{\varphi}$ is linear, thus any harmonic map is
biharmonic. We call {\it proper biharmonic} the non-harmonic
biharmonic maps.

\indent In this paper we shall focus our attention on {\it biharmonic
submanifolds}, i.e.  on submanifolds such that the inclusion map is
a biharmonic map. In this context, a proper biharmonic submanifold is
a non-minimal biharmonic submanifold.

The first ambient spaces to look for proper biharmonic
submanifolds are the spaces of constant sectional curvature $c$,
which we shall denote by $\e^n(c)$, and the first class of
submanifolds to be studied is that of the hypersurfaces. The full
classification of proper biharmonic hypersurfaces in  $\e^n(c)$,
for any $n\geq 3$, is not known yet  and, up to now, these are the
results obtained:
\begin{itemize}
  \item biharmonic hypersurfaces in $\r^n$, $n=3,4$, are  minimal
  \cite{BYC_ISH,THTV};
  \item biharmonic surfaces in $\h^3$ are minimal \cite{RCSMCO2};
  \item biharmonic surfaces in $\s^3$ are open parts
  of the hypersphere $\s^2(\rrm)$ \cite{RCSMCO1}.
\end{itemize}

The aim of this paper is to go further with the classification of
compact proper biharmonic hypersurfaces in $\e^n(c)$. This study
will conduce to the following classification of proper biharmonic
compact hypersurfaces in $\s^4$.

\begin{theorem*}
The only proper biharmonic compact hypersurfaces of $\s^4$ are the
hypersphere $\s^3(\rrm)$ and the torus
$\s^1(\rrm)\times\s^2(\rrm)$.
\end{theorem*}

The strategy to prove the theorem consists in proving that
 proper biharmonic  hypersurfaces in $4$-dimensional space forms
 have constant mean
curvature. This is done by dividing the study according to the
number of distinct principal curvatures.

The simplest assumption that $M$ is an umbilical hypersurface,
i.e. all principal curvatures are equal, gives an immediate picture.
 In fact, if $M$ is a proper biharmonic umbilical
hypersurface in $\s^{m+1}$, then it is an open part of
$\s^m(\rrm)$. Moreover, there exist no proper biharmonic umbilical
hypersurfaces in $\r^{m+1}$ or in the hyperbolic space $\h^{m+1}$
(see \cite{RCSMCO2}).

For biharmonic hypersurfaces with at most two distinct principal
curvatures the property of having constant mean curvature was
proved in \cite{ID}, for the Euclidean case, and in \cite{ABSMCO}
for any space form. This property proved to be the main ingredient
for the following complete classification of proper biharmonic
hypersurfaces with at most two distinct principal curvatures in
the Euclidean sphere.
\begin{theorem}[\cite{ABSMCO}]\label{th: classif_hypersurf_2_curv_princ}
Let $M^m$ be a proper biharmonic hypersurface with at most two
distinct principal curvatures in $\s^{m+1}$. Then $M$ is an open
part of $\s^{m}(\rrm)$ or of $\s^{m_1}(\rrm)\times
\s^{m_2}(\rrm)$, $m_1+m_2=m$, $m_1\neq m_2$.
\end{theorem}

In this paper we first prove that there exist no compact proper
 biharmonic hypersurfaces, of constant
mean curvature, with three distinct principal curvatures in $\s^n$
 (Theorem~\ref{th: non-exist_isoparam_3princ}).
Secondly, we show that biharmonic hypersurfaces of $\e^4(c)$ must
have constant mean curvature (Theorem~\ref{th:
curb_med_const_3_curv_princ}).

These two results, together with Theorem~\ref{th:
classif_hypersurf_2_curv_princ}, give, as a consequence, the main
result of the paper.

For an up-to-date bibliography on biharmonic maps we refer the reader to
\cite{BIB}.

\section{Preliminaries}

Let $\varphi:M\to\e^n(c)$  be the canonical inclusion of a
submanifold $M$ in a constant sectional curvature $c$ manifold,
$\e^n(c)$. The expressions assumed by the tension and bitension
fields are
\begin{equation}\label{eq: tens+bitens E^n(c)}
\tau(\varphi)=mH,\qquad\qquad \tau_2(\varphi)=-m(\Delta H-mcH),
\end{equation}
where $H$ denotes the mean curvature vector field of $M$ in
$\e^n(c)$, while $\Delta$ is the rough Laplacian on
$\varphi^{-1}T\e^n(c)$.

The following characterization result, obtained in \cite{BYC2} and
\cite{CO1} by splitting the bitension field in its normal and
tangent components, represents the main tool in the study of
proper biharmonic submanifolds in space forms.

\begin{theorem}[\cite{BYC2, CO1}]\label{th: bih subm E^n(c)}
The canonical inclusion $\varphi:M^m\to\e^n(c)$ of a submanifold $M$
in an $n$-dimensional space form $\e^n(c)$ is biharmonic if and
only if
\begin{equation}\label{caract_bih_spheres}
\left\{
\begin{array}{l}
\ -\Delta^\perp H-\trace B(\cdot,A_H\cdot)+mcH=0,
\\ \mbox{} \\
\ 2\trace A_{\nabla^\perp_{(\cdot)}H}(\cdot)
+\frac{m}{2}\grad(\vert H \vert^2)=0,
\end{array}
\right.
\end{equation}
where $A$ denotes the Weingarten operator, $B$ the second
fundamental form, $H$ the mean curvature vector field,
$\nabla^\perp$ and $\Delta^\perp$ the connection and the Laplacian
in the normal bundle of $M$ in $\e^n(c)$.

Moreover, if $M$ is a hypersurface of $\e^{m+1}(c)$, then $M$ is
proper biharmonic if and only if
\begin{equation}\label{caract_bih_hipersurf_spheres}
\left\{
\begin{array}{l}
\Delta^\perp H-(mc-|A|^2)H=0,
\\ \mbox{} \\
\ 2A\big(\grad (|H|)\big)+m|H|\grad(|H|)=0.
\end{array}
\right.
\end{equation}
\end{theorem}

The generalized Clifford torus
$\s^{m_1}(\rrm)\times\s^{m_2}(\rrm)$, $m_1+m_2=m$, $m_1\neq m_2$,
was the first example of proper biharmonic hypersurface in
$\s^{m+1}$ (see \cite{GYJ1}). Then, in \cite{RCSMCO1} it was
proved that the only proper biharmonic hypersphere $\s^{m}(a)$,
$a\in(0,1)$, in $\s^{m+1}$ is $\s^{m}(\rrm)$.

Inspired by these fundamental examples, in \cite{RCSMCO2}, the
authors presented two methods for constructing proper biharmonic
submanifolds of codimension greater  than $1$ in $\s^n$.

\begin{theorem}[Composition property, \cite{RCSMCO2}] \label{th: rm_minim}
Let $M$ be a minimal submanifold of\, $\s^{n-1}(a)\subset\s^n$.
Then $M$ is proper biharmonic in $\s^{n}$ if and only if $a=\rrm$.
\end{theorem}

\begin{theorem}[Product composition property, \cite{RCSMCO2}]\label{th:hipertor}
Let $M_1^{m_1}$ and $M_2^{m_2}$ be two minimal submanifolds of
$\s^{n_1}(r_1)$ and $\s^{n_2}(r_2)$, respectively, where
$n_1+n_2=n-1$, $r_1^2+r_2^2=1$. Then $M_1\times M_2$ is proper
biharmonic in $\s^n$ if and only if $r_1=r_2=\rrm$ and $m_1\neq
m_2$.
\end{theorem}

\subsection{Isoparametric hypersurfaces}

We recall that a hypersurface $M^m$ in $\s^{m+1}$ is said to be
{\it isoparametric} of type $\ell$ if it has constant principal
curvatures $k_1 > \ldots
> k_\ell$ with respective constant multiplicities $m_1, \ldots
,m_\ell$, $m = m_1 + m_2 + \ldots + m_\ell$. It is known that the
number $\ell$ is either $1, 2, 3, 4$ or $6$.  For $\ell \leq 3$ we
have the following classification of compact isoparametric
hypersurfaces.

If $\ell  = 1$, then $M$ is totally umbilical.

If $\ell = 2$, then $M = \s^{m_1}(r_1)\times \s^{m_2}(r_2)$,
$r_1^2 + r_2^2 = 1$ (see \cite{RY}).

If $\ell = 3$, then $m_1 = m_2 = m_3 = 2^q$, $q=0,1,2,3$ (see
\cite{CAR}).

Moreover, there exists an angle $\theta$, $0 < \theta <
\frac{\pi}{\ell}$ , such that
\begin{equation}\label{eq-kalpha}
k_\alpha = \cot\big(\theta + \frac{(\alpha-1)\pi}{\ell}\big),\quad
\alpha = 1,\ldots, \ell.
\end{equation}

In the next section we shall need the following results on
isoparametric hypersurfaces

\begin{theorem}[\cite{CHA}]\label{th: CMC_SMC_3dist_princ_curv}
A compact hypersurface $M^m$ of constant scalar curvature $s$ and
constant mean curvature $|H|$ in $\s^{m+1}$ is isoparametric
provided it has $3$ distinct principal curvatures everywhere.
\end{theorem}

\begin{theorem}[\cite{CHA1}]\label{th: Chang_S4}
Any compact hypersurface with constant scalar and mean curvature
in $\s^4$ is isoparametric.
\end{theorem}

We end this section by recalling
\begin{proposition}[\cite{ABSMCO}]\label{th: curb_scal_hyp}
Let $M^m$ be a proper biharmonic hypersurface with constant mean
curvature $|H|^2=k$ in $\s^{m+1}$. Then $M$ has constant scalar
curvature,
$$
s=m^2(1+k)-2m.
$$
\end{proposition}

\section{Biharmonic hypersurfaces with three distinct principal curvatures in spheres}

Using the classification result on isoparametric hypersurfaces we
can prove the following non-existence result for biharmonic
hypersurfaces with $3$ distinct principal curvatures.

\begin{theorem}\label{th: non-exist_isoparam_3princ}
There exist no compact proper biharmonic hypersurfaces of constant
mean curvature and with three distinct principal curvatures in the
unit Euclidean sphere.
\end{theorem}
\begin{proof} By Proposition \ref{th: curb_scal_hyp}, a proper biharmonic hypersurface $M$
with constant mean curvature $|H|^2=k$ in $\s^{m+1}$ has constant
scalar curvature. Since $M$ is compact with $3$ distinct principal
curvatures and has constant mean and scalar curvature, from
Theorem~\ref{th: CMC_SMC_3dist_princ_curv}, it results that $M$ is
isoparametric with $\ell=3$ in $\s^{m+1}$. Now, taking into
account \eqref{eq-kalpha}, there exists $\theta\in(0,\pi/3)$ such
that
$$
k_1=\cot \theta,\qquad k_2 = \cot\big(\theta +
\frac{\pi}{3}\big)=\frac{k_1-\sqrt 3}{1+\sqrt{3}k_1},\qquad k_3 =
\cot\big(\theta + \frac{2\pi}{3}\big)=\frac{k_1+\sqrt
3}{1-\sqrt{3}k_1}.
$$
Thus, from Cartan's classification, the square of the norm of the
shape operator is
\begin{equation}\label{eq-akalpha}
|A|^2=2^q(k_1^2+k_2^2+k_3^2)
=2^q\frac{9k_1^6+45k_1^2+6}{(1-3k_1^2)^2}
\end{equation}
and $m=3\cdot2^q$, $q=0,1,2,3$. On the other hand, since $M$ is
biharmonic of constant mean curvature, from
\eqref{caract_bih_hipersurf_spheres},
$$
|A|^2=m=3\cdot2^q.
$$
The last equation, together with \eqref{eq-akalpha}, implies that
$k_1$ is a solution of $3k_1^6-9k_1^4+21k_1^2+1=0$, which is an
equation with no real roots.

\end{proof}

\begin{remark}
Investigating the biharmonicity of compact isoparametric
hypersurfaces in Euclidean spheres, in \cite{IIU} the authors
proved that the only compact proper biharmonic isoparametric
hypersurfaces of $\s^{m+1}$ are the hypersphere $\s^m(\rrm)$ and
the generalized Clifford torus
$\s^{m_1}(\rrm)\times\s^{m_2}(\rrm)$, $m_1+m_2=m$, $m_1\neq m_2$.
\end{remark}

We shall now concentrate on biharmonic hypersurfaces in $\e^4(c)$.
Using B-Y.~Chen techniques (see also \cite{DKP}) we prove

\begin{theorem}\label{th: curb_med_const_3_curv_princ}
Let $M^3$ be a biharmonic hypersurface of the space form
$\e^4(c)$. Then $M$ has constant mean curvature.
\end{theorem}
\begin{proof}

Suppose that $|H|$ is not constant on $M$. Then there exists an
open subset $U$ of $M$ such that $\grad_p|H|^2\neq 0$, for all
$p\in U$. By eventually restraining $U$ we can suppose that
$|H|>0$ on $U$, and thus $\grad_p|H|\neq 0$, for all $p\in U$. If
$U$ has at most two distinct principal curvatures, then , by Theorem~4.1 in \cite{ABSMCO},
 we conclude that its mean
curvature is constant and we have a contradiction. Then there
exists a point in $U$ with three distinct principal curvatures.
This implies the existence of an open neighborhood of points with
three distinct principal curvatures and we can suppose, by
restraining $U$, that all its points have three distinct principal
curvatures. On $U$ we can consider the unit section in the normal
bundle $\eta=\frac{H}{|H|}$ and denote by $f=|H|$ the mean
curvature function of $U$ in $\e^{m+1}(c)$ and by $k_i$,
$i=1,2,3$, its principal curvatures w.r.t. $\eta$.

Conclusively, the hypothesis for $M$ to be proper biharmonic with
at most three distinct principal curvatures in $\e^{m+1}(c)$ and
non-constant mean curvature, implies the existence of an open
connected subset $U$ of $M$, with
\begin{equation}\label{cond_f_3}
\left\{
\begin{array}{l}
\grad_p f\neq 0,
\\ \mbox{} \\
f(p)>0,
\\ \mbox{} \\
k_1(p)\neq k_2(p)\neq k_3(p)\neq k_1(p), \qquad\forall p\in U.
\end{array}
\right.
\end{equation}

We shall contradict the condition $\grad_p f\neq 0$, for all $p\in
U$.

Since $M$ is proper biharmonic in $\e^{4}(c)$, from
\eqref{caract_bih_hipersurf_spheres} we have
\begin{equation}\label{eq: caract_bih_spheres_hypersurf}
\left\{
\begin{array}{l}
\ \Delta f=(3c-|A|^2)f,
\\ \mbox{} \\
A(\grad f)=-\frac{3}{2}f\grad f.
\end{array}
\right.
\end{equation}

Consider now $X_1=\displaystyle{\frac{\grad f}{|\grad f|}}$ on
$U$. Then $X_1$ is a principal direction with principal curvature
$k_1=-\frac{3}{2}f$. Recall that $3f=k_1+k_2+k_3$, thus
\begin{equation}\label{eq: sum_curv23}
k_2+k_3=\frac{9}{2}f.
\end{equation}

We shall use the moving frames method and denote by $X_1, X_2,
X_3$ the orthonormal frame field of principal directions and by
$\{\omega^a\}_{a=1}^3$ the dual frame field of $\{X_a\}_{a=1}^3$
on $U$.

Obviously,
\begin{equation}\label{eq: X_alpha(f)=0}
X_i(f)=\langle X_i, \grad f\rangle=|\grad f|\langle
X_i,X_1\rangle=0, \qquad i=2, 3,
\end{equation}
 thus
\begin{equation}\label{eq: grad_f}
\grad f=X_1(f)X_1.
\end{equation}

We write
$$
\nabla X_a=\omega_a^bX_b, \qquad \omega_a^b\in C(T^*U).
$$
From the Codazzi equations for $M$, for distinct $a, b, d=1, 2, 3$,  we get
\begin{eqnarray}\label{eq: Codazzi_hyper1_3c}
X_a(k_b)=(k_a-k_b)\omega_a^b(X_b)
\end{eqnarray}
and
\begin{eqnarray}\label{eq: Codazzi_hyper2_3c}
(k_b-k_d)\omega^d_b(X_a)=(k_a-k_d)\omega^d_a(X_b).
\end{eqnarray}

Consider now in \eqref{eq: Codazzi_hyper1_3c}, $a=1$ and $b=i$
and, respectively, $a=i$ and $b=j$ with $i\neq j$. We obtain
\begin{eqnarray*}\label{eq: omega_fund_3c}
\omega^1_i(X_i)&=&\frac{X_1(k_i)}{k_i-k_1}
\end{eqnarray*}
and
\begin{eqnarray*}\label{eq: omegaij(xj)}
\omega^i_j(X_j)&=&\frac{X_i(k_j)}{k_j-k_i}.
\end{eqnarray*}

For $a=i$ and $b=1$, as $X_i(k_1)=0$, \eqref{eq:
Codazzi_hyper1_3c} leads to  $\omega^1_i(X_1)=0$ and we can write
\begin{equation*}\label{eq: omega1a(x1)_3c}
\omega^1_a(X_1)=0,\quad a=1, 2, 3.
\end{equation*}

Notice that, since $X_i(f)=0$, then $\langle[X_i,X_j],
X_1\rangle=0$, thus $\omega_1^j(X_i)=\omega_1^i(X_j)$. Now, from
\eqref{eq: Codazzi_hyper2_3c}, for $a=1$, $b=i$ and $d=j$, with
$i\neq j$,
 we get
\begin{equation*}\label{eq: omega1b(xa)_3c}
\omega^1_2(X_3)=\omega^2_3(X_1)=\omega^3_1(X_2)=0.
\end{equation*}

The structure $1$-forms are thus determined by the following set
of relations

\begin{equation}\label{eq: 1-form_conn}
\left\{
\begin{array}{lll}
\omega^1_2(X_1)=0, &
\omega^1_2(X_2)=\frac{X_1(k_2)}{k_2+\frac{3}{2}f}=\alpha_2, &
\omega^1_2(X_3)=0,
\\ \mbox{} \\
\omega^1_3(X_1)=0, & \omega^1_3(X_2)=0, &
\omega^1_3(X_3)=\frac{X_1(k_3)}{k_3+\frac{3}{2}f}=\alpha_3,
\\ \mbox{} \\
\omega^2_3(X_1)=0,&
\omega^2_3(X_2)=\frac{X_3(k_2)}{k_3-k_2}=\beta_2, &
\omega^2_3(X_3)=\frac{X_2(k_3)}{k_3-k_2}=\beta_3,
\end{array}
\right.
\end{equation}

In order to express the first condition in \eqref{eq:
caract_bih_spheres_hypersurf}, by using \eqref{eq: sum_curv23}, we
compute
\begin{eqnarray}\label{eq: norm_A2}
|A|^2&=&k_1^2+k_2^2+k_3^2\nonumber\\
&=& k_1^2+(k_2+k_3)^2-2k_2k_3\\
&=&\frac{45}{2}f^2-2K,\nonumber
\end{eqnarray}
where $K$ denotes the product $k_2k_3$. From \eqref{eq: grad_f} we
deduce that
\begin{eqnarray}\label{eq: Delta_f}
\Delta f&=&-\Div(\grad f)=-\Div(X_1(f)X_1)=-X_1(X_1(f))-X_1(f)\Div
X_1\nonumber
\\&=&-X_1(X_1(f))+X_1(f)(\omega^1_2(X_2)+\omega^1_3(X_3))\nonumber\\
&=&-X_1(X_1(f))+X_1(f)(\alpha_2+\alpha_3).
\end{eqnarray}
Now, by using \eqref{eq: norm_A2} and \eqref{eq: Delta_f}, the
equation $\Delta f=(3c-|A|^2)f$ becomes
\begin{equation}\label{eq: II}
X_1(X_1(f))-X_1(f)(\alpha_2+\alpha_3)+(2K+3c-\frac{45}{2}f^2)f=0.
\end{equation}
We also compute
\begin{eqnarray}\label{eq: comm 12_13}
[X_1,X_i]&=&\nabla_{X_1}X_i-\nabla_{X_i}X_1=\langle
\nabla_{X_1}X_i,X_1\rangle X_1-\langle \nabla_{X_i}X_1,X_i\rangle
X_i\nonumber\\
&=&\omega^1_i(X_i)X_i=\alpha_i X_i.
\end{eqnarray}

We shall now use the Gauss equation
\begin{eqnarray}\label{eq: Gauss_EQ}
\langle R^{\e^4(c)}(X,Y)Z,W\rangle&=&\langle
R(X,Y)Z,W\rangle\nonumber\\&&+\langle B(X,Z),B(Y,W)\rangle-\langle
B(X,W),B(Y,Z)\rangle.
\end{eqnarray}
From \eqref{eq: Gauss_EQ} we have:
\begin{itemize}
\item for $X=W=X_1$ and $Y=Z=X_i$
\begin{equation}\label{eq: Gauss1}
\left\{
\begin{array}{l}
X_1(\alpha_2)=\alpha_2^2+c-\frac{3}{2}fk_2,
\\\mbox{}\\
X_1(\alpha_3)=\alpha_3^2+c-\frac{3}{2}fk_3;
\end{array}
\right.
\end{equation}
\item  for $X=W=X_2$ and $Y=Z=X_3$
\begin{equation}\label{eq: Gauss2}
K+c=X_2(\beta_3)-X_3(\beta_2)
-\alpha_2\alpha_3-\beta_2^2-\beta_3^2;
\end{equation}
\item for $Y=W=X_3$, $X=X_2$ and $Z=X_1$ and,
respectively, for $X=W=X_2$, $Y=X_3$ and $Z=X_1$
\begin{equation}\label{eq: Gauss3}
\left\{
\begin{array}{l}
X_2(\alpha_3)=\beta_3(\alpha_3-\alpha_2),
\\\mbox{}\\
X_3(\alpha_2)=\beta_2(\alpha_3-\alpha_2);
\end{array}
\right.
\end{equation}
\item  for $X=W=X_2$, $Y=X_1$ and $Z=X_3$,
and, respectively, for $X=W=X_3$, $Y=X_2$ and $Z=X_1$
\begin{equation}\label{eq: Gauss4}
\left\{
\begin{array}{l}
X_1(\beta_2)=\alpha_2\beta_2,
\\\mbox{}\\
X_1(\beta_3)=\alpha_3\beta_3.
\end{array}
\right.
\end{equation}
\end{itemize}
Notice now that, from \eqref{eq: X_alpha(f)=0} and \eqref{eq: comm
12_13},
\begin{equation}\label{eq: X_iX_1f}
X_i(X_1(f))=-[X_1,X_i]f+X_1(X_i(f))=-\alpha_i X_i(f)+X_1(X_i(f))=0
\end{equation}
and
\begin{equation}\label{eq: X_iX_1X_1f}
X_i(X_1(X_1(f)))=0.
\end{equation}
Also, since $K=\frac{(k_2+k_3)^2-(k_3-k_2)^2}{4}$ we obtain
\begin{equation}\label{eq: X_iK}
\left\{
\begin{array}{l}
X_2(K)=-(k_3-k_2)^2\beta_3,
\\\mbox{}\\
X_3(K)=(k_3-k_2)^2\beta_2.
\end{array}
\right.
\end{equation}
We differentiate \eqref{eq: II} along $X_2$ and $X_3$ and use
\eqref{eq: Gauss3}, \eqref{eq: X_iX_1f}, \eqref{eq: X_iX_1X_1f}
and \eqref{eq: X_iK}. We get
\begin{equation}\label{eq: diff_consecII_12}
\left\{
\begin{array}{l}
X_2(\alpha_2)=-\beta_3(\alpha_3-\alpha_2)-\frac{2f}{X_1(f)}(k_3-k_2)^2\beta_3,
\\\mbox{}\\
X_3(\alpha_3)=-\beta_2(\alpha_3-\alpha_2)+\frac{2f}{X_1(f)}(k_3-k_2)^2\beta_2.
\end{array}
\right.
\end{equation}

We intend to prove that $X_i(k_j)=0$, $i,j=2,3$. In order to do
this we apply $[X_1,X_2]=\alpha_2 X_2$ to the quantity $\alpha_2$.
On one hand, from \eqref{eq: diff_consecII_12}, we get
\begin{equation}\label{eq: beta3=0_1}
[X_1,X_2]\alpha_2=\alpha_2X_2(\alpha_2)=
\beta_3\Big\{-\alpha_2\alpha_3+\alpha_2^2-\frac{2f}{X_1(f)}(k_3-k_2)^2\alpha_2\Big\}.
\end{equation}
On the other hand, by using \eqref{eq: Gauss1} and \eqref{eq:
diff_consecII_12}, we obtain
\begin{eqnarray}\label{eq: beta3=0_2}
[X_1,X_2]\alpha_2&=&X_1(X_2(\alpha_2))-X_2(X_1(\alpha_2))\nonumber\\
&=&\beta_3\Big\{-2\alpha_3^2-\alpha_2^2+3\alpha_2\alpha_3+
\frac{2f}{X_1(f)}\big[-2(k_3-k_2)X_1(k_3-k_2)\nonumber\\
&&+(k_3-k_2)^2(2\alpha_2-\alpha_3)\big]-2X_1\Big(\frac{f}{X_1(f)}\Big)(k_3-k_2)^2\Big\}.
\end{eqnarray}
By putting together \eqref{eq: beta3=0_1} and \eqref{eq:
beta3=0_2} we either have $\beta_3=0$ or
\begin{eqnarray}\label{eq: beta3=0_3}
X_1\Big(\frac{f}{X_1(f)}\Big)=-\frac{(\alpha_3-\alpha_2)^2}{(k_3-k_2)^2}
+\frac{f}{X_1(f)}\Big(3\alpha_2-\alpha_3-2\frac{X_1(k_3-k_2)}{k_3-k_2}\Big).
\end{eqnarray}
Moreover,
\begin{eqnarray*}
X_2\Big(\frac{X_1(k_3-k_2)}{k_3-k_2}\Big)&=&
-\frac{1}{k_3-k_2}[X_1,X_2](k_3-k_2)+X_1\Big(\frac{X_2(k_3-k_2)}{k_3-k_2}\Big)\nonumber\\
&=&2(\alpha_3-\alpha_2)\beta_3.
\end{eqnarray*}
Suppose that $\beta_3\neq 0$, differentiate \eqref{eq: beta3=0_3}
along $X_2$ and use \eqref{eq: Gauss3} and \eqref{eq:
diff_consecII_12}. We get
\begin{equation}\label{eq: beta3=0_4}
2(\alpha_3-\alpha_2)=-\frac{f}{X_1(f)}(k_3-k_2)^2.
\end{equation}
We differentiate now \eqref{eq: beta3=0_4} along $X_2$ and obtain
\begin{equation}\label{eq: beta3=0_5}
\alpha_3-\alpha_2=-2\frac{f}{X_1(f)}(k_3-k_2)^2,
\end{equation}
and since $k_2\neq k_3$ the equations \eqref{eq: beta3=0_4} and
\eqref{eq: beta3=0_5} lead to a contradiction.

Analogously, by using the symmetry of the equations in $X_2$ and
$X_3$, we immediately prove that $\beta_2=0$.

We rewrite equations \eqref{eq: Gauss1} in the form
\begin{equation}\label{eq: Gauss1_reloaded}
\left\{
\begin{array}{l}
X_1(X_1(k_2))=\frac{21}{2}\alpha_2
X_1(f)+2(K+c)(k_3+\frac{3}{2}f)+(c-\frac{3}{2}fk_2)(k_2+\frac{3}{2}f),
\\\mbox{}\\
X_1(X_1(k_3))=\frac{21}{2}\alpha_3
X_1(f)+2(K+c)(k_2+\frac{3}{2}f)+(c-\frac{3}{2}fk_3)(k_3+\frac{3}{2}f),
\end{array}
\right.
\end{equation}
and by summing up we obtain
\begin{equation}\label{eq: f''_from_Gauss1}
X_1(X_1(f))=\frac{7}{3}X_1(f)(\alpha_2+\alpha_3)+f(4K+5c-9f^2).
\end{equation}
Now, by using \eqref{eq: II} and \eqref{eq: f''_from_Gauss1} we
obtain
\begin{equation}\label{eq: X_1f_1}
X_1(f)(\alpha_2+\alpha_3)=f\big(-\frac{9}{2}K-6c+\frac{189}{8}f^2\big).
\end{equation}
We replace \eqref{eq: X_1f_1} in \eqref{eq: f''_from_Gauss1} and
get
\begin{equation}\label{eq: X_1X_1f_1}
X_1(X_1(f))=f\big(-\frac{13}{2}K-9c+\frac{369}{8}f^2\big).
\end{equation}

In order to get another relation on $f$ and $K$ we first use
\eqref{eq: Gauss2}, \eqref{eq: Gauss1}, \eqref{eq: sum_curv23},
\eqref{eq: 1-form_conn} and determine
\begin{eqnarray}\label{eq: X_1K}
  X_1(K) &=& -X_1(\alpha_2\alpha_3) \\
  &=& -(\alpha_2\alpha_3+c)(\alpha_2+\alpha_3)+\frac{3}{2}f(\alpha_2k_3+\alpha_3k_2) \nonumber\\
  &=& (K+9f^2)(\alpha_2+\alpha_3)-\frac{27}{4}fX_1(f).\nonumber
\end{eqnarray}

By differentiating \eqref{eq: X_1f_1} along $X_1$, and by
using \eqref{eq: X_1K}, \eqref{eq: X_1X_1f_1}, \eqref{eq: Gauss1},
\eqref{eq: X_1f_1} we get
\begin{equation}\label{eq: X_1f_2}
X_1(f)\Big(\frac{13}{2}K+10c-108f^2\Big)=f(\alpha_2+\alpha_3)\Big(\frac{13}{2}K+15c-\frac{441}{4}f^2\Big).
\end{equation}
We multiply \eqref{eq: X_1f_2} first by $X_1(f)$ and secondly by
$\alpha_2+\alpha_3$ and, by using \eqref{eq: X_1f_1}, we get
{\small
\begin{equation}\label{eq: (X_1f)2_(a2a3)^2} \left\{
\begin{array}{l}
(X_1(f))^2(\frac{13}{2}K+10c-108f^2)=f^2(-\frac{9}{2}K-6c+\frac{189}{8}f^2)(\frac{13}{2}K+15c-\frac{441}{4}f^2),
\\\mbox{}\\
(\frac{13}{2}K+10c-108f^2)(-\frac{9}{2}K-6c+\frac{189}{8}f^2)=(\alpha_2+\alpha_3)^2(\frac{13}{2}K+15c-\frac{441}{4}f^2).
\end{array}
\right.
\end{equation}}

Differentiating \eqref{eq: X_1f_1} along $X_1$, and using
\eqref{eq: X_1K}, \eqref{eq: X_1X_1f_1}, \eqref{eq: Gauss1},
\eqref{eq: X_1f_1}, \eqref{eq: (X_1f)2_(a2a3)^2}, we obtain
\begin{align}\label{eq: first_pol}
27f^2(4044800c^3 - 49579440c^2 f^2 + 187840944 c f^4 -
          254205945 f^6)\nonumber\\ -
    6 (51200 c^3 - 19600320 c^2 f^2 + 119328660 c f^4 -
          80969301 f^6) K\\ -
    208 (2240 c^2 - 108396 c f^2 - 285363 f^4) K^2\nonumber\\ +
    2704 (16 c - 2277 f^2)K^3\nonumber\\ + 140608 K^4=0.\nonumber
\end{align}

Consider now $\gamma=\gamma(t)$, $t\in I$, to be an integral curve
of $X_1$ passing through $p=\gamma(t_0)$. Since $X_2(f)=X_3(f)=0$
and $X_2(K)=X_3(K)=0$ and $X_1(f)\neq 0$, we can write $t=t(f)$ in
a neighborhood of $f_0=f(t_0)$ and thus consider $K=K(f)$.

Notice that if $\frac{13}{2}K+15c-\frac{441}{4}f^2=0$ or
$10c+\frac{13}{2}K-108f^2=0$, then from \eqref{eq: first_pol} the
function $f$ results to be the solution of a polynomial equation
of eighth degree with constant coefficients and we would get to a
contradiction. Thus, from \eqref{eq: (X_1f)2_(a2a3)^2} we have
that
\begin{equation}\label{eq: (X_1f)2_(a2a3)^2...2}
\left\{
\begin{array}{l}
\big(\frac{df}{dt}\big)^2=\frac{f^2(-\frac{9}{2}K-6c+\frac{189}{8}f^2)(\frac{13}{2}K+15c-\frac{441}{4}f^2)}{\frac{13}{2}K+10c-108f^2},
\\\mbox{}\\
(\alpha_2+\alpha_3)^2=\frac{(\frac{13}{2}K+10c-108f^2)(-\frac{9}{2}K-6c+\frac{189}{8}f^2)}{\frac{13}{2}K+15c-\frac{441}{4}f^2}.
\end{array}
\right.
\end{equation}

We can now compute $\frac{dK}{df}$ by using \eqref{eq:
(X_1f)2_(a2a3)^2...2}, \eqref{eq: X_1K} and \eqref{eq: X_1f_1},

\begin{eqnarray}\label{eq: dK/df}
\frac{dK}{df}&=&\frac{dK}{dt}\frac{dt}{df}\nonumber\\
&=&\frac{(K+9f^2)\frac{df}{dt}(\alpha_2+\alpha_3)}{(\frac{df}{dt})^2}-\frac{27}{4}f\nonumber\\
&=&\frac{(K+9f^2)(\frac{13}{2}K+10c-108f^2)}{f(\frac{13}{2}K+15c-\frac{441}{4}f^2)}-\frac{27}{4}f.
\end{eqnarray}

The next step consists in differentiating \eqref{eq: first_pol}
with respect to $f$. By substituting $\frac{dK}{df}$ from
\eqref{eq: dK/df} we get another polynomial equation in $f$ and
$K$, of fifth degree in $K$. We eliminate $K^5$ between this new
polynomial equation and \eqref{eq: first_pol}. The result
constitutes a polynomial equation in $f$ and $K$, of fourth degree
in $K$. In a similar way, by using \eqref{eq: first_pol} and its
consequences we are able to gradually eliminate $K^4$, $K^3$,
$K^2$ and $K$ and we are led to a polynomial equation with
constant coefficients in $f$. Thus $f$ results to be constant and
we conclude.
\end{proof}

In \cite{THTV} the authors proved that there exist no proper
biharmonic hypersurfaces in $\r^4$. By using Theorem \ref{th:
curb_med_const_3_curv_princ} we reobtain this result and we can extend
it to the $4$-dimensional hyperbolic space.
\begin{theorem}
There exist no proper biharmonic hypersurfaces in $\r^4$ or
$\h^{4}$.
\end{theorem}

\begin{proof}
Suppose that $M^3$ is a proper biharmonic hypersurface in $\r^{4}$
or $\h^{4}$. From Theorem \ref{th: curb_med_const_3_curv_princ},
the mean curvature of $M$ is constant, and applying Theorem
\ref{th: bih subm E^n(c)} we obtain $|A|^2=0$ or $|A|^2=-3$,
respectively, and we conclude.
\end{proof}

We can now state the main result of the paper.

\begin{theorem}\label{th: hyper_S4}
The only proper biharmonic compact hypersurfaces of $\s^4$ are the
hypersphere $\s^3(\rrm)$ and the torus
$\s^1(\rrm)\times\s^2(\rrm)$.
\end{theorem}

\begin{proof}
Suppose that $M^3$ is a compact proper biharmonic hypersurface of
$\s^4$. From Theorem \ref{th: curb_med_const_3_curv_princ} it
results that $M$ has constant mean curvature and, since it
satisfies the hypotheses of Proposition \ref{th: curb_scal_hyp},
we conclude that it also has constant scalar curvature. We can
thus apply Theorem \ref{th: Chang_S4} and it results that $M$ is
isoparametric in $\s^4$. From Theorem \ref{th:
non-exist_isoparam_3princ} we get that $M$ cannot be isoparametric
with $\ell=3$, and by using Theorem \ref{th:
classif_hypersurf_2_curv_princ} we conclude the proof.

\end{proof}

\begin{remark}
Since we have achieved the classification of proper biharmonic
curves (see \cite{RCSMCO2}) and compact proper biharmonic
hypersurfaces, our interest is now in
 proper biharmonic surfaces of $\s^4$. Using Theorem~\ref{th: rm_minim} and a well known
 result of Lawson \cite{HBL},
  it was proved, in \cite{RCSMCO2}, the
existence of closed orientable embedded proper biharmonic surfaces
of arbitrary genus in $\s^4$. All this surfaces are  minimal surfaces of $\s^3(\rrm)$.
Moreover, Theorem~\ref{th:hipertor} cannot be applied due to the dimensions
involved, thus it does not generate examples in this case. Then it
is natural to propose the following

\subsection*{Open problem} {\it Are there other proper biharmonic surfaces
in $\s^4$, apart from the minimal surfaces of $\s^3(\rrm)$?}

\end{remark}

\end{document}